%% file: cpt_opt.tex
\documentclass[12pt]{article}
\usepackage{fullpage,graphicx,psfrag,amsmath,amsfonts,verbatim}
\usepackage[small,bf]{caption}\usepackage{amssymb}
\usepackage{fullpage}
\usepackage{enumitem}
\usepackage{tikz}
\usepackage{listings}
\usepackage{authblk}
\usepackage{textcomp}

\graphicspath{{figures/}}

\DeclareFixedFont{\ttb}{T1}{txtt}{bx}{n}{12} 
\DeclareFixedFont{\ttm}{T1}{txtt}{m}{n}{12}  

\usepackage{color}
\definecolor{deepblue}{rgb}{0,0,0.5}
\definecolor{deepgreen}{rgb}{0,0.5,0}

\usepackage{hyperref}
\usepackage{wasysym}
\hypersetup{
colorlinks   = true,    
urlcolor     = blue,    
linkcolor    = blue,    
citecolor    = red      
}

\input defs.tex
\title{Portfolio Optimization with Cumulative Prospect Theory Utility
via Convex Optimization}

\author[1]{Eric Luxenberg\footnote{Equal contribution.}}
\newcommand\CoAuthorMark{\footnotemark[\arabic{footnote}]} 

\author[2]{Philipp Schiele\protect\CoAuthorMark}
\author[1]{Stephen Boyd}
\affil[1]{Department of Electrical Engineering, Stanford University}
\affil[2]{Department of Statistics, Ludwig-Maximilians-Universität München}

\begin{document}
\maketitle

\begin{abstract}
We consider the problem of choosing a portfolio that maximizes the
cumulative prospect theory (CPT) utility on an empirical
distribution of asset returns.
We show that while CPT utility is not a concave function of the portfolio
weights, it can be expressed as a difference of two functions.
The first term is the composition of a convex function with concave arguments
and the second term a composition of a convex function with convex arguments.
This structure allows us to derive a global lower bound, or minorant, on the
CPT utility, which we can use in a minorization-maximization (MM) algorithm for
maximizing CPT utility.
We further show that the problem is amenable to a simple convex-concave (CC)
procedure which iteratively maximizes a local approximation.
Both of these methods can handle small and medium size problems, and 
complex (but convex) portfolio constraints.
We also describe a simpler method that scales to larger problems, but
handles only simple portfolio constraints.
\end{abstract}

\newpage
\tableofcontents
\newpage

\section{Introduction}
\subsection{Cumulative prospect theory}
Analysis of decision-making under uncertainty has long been dominated by 
von Neumann-Morgenstern (VNM) utility maximization~\cite{utility},
which takes rational behavior as a fundamental assumption.
Kahneman and Tversky~\cite{kahneman1979prospect} observed that the
VNM theory fails to explain actual human decision-making behavior in
many settings. The subsequently introduced prospect theory (PT) formalizes loss
aversion and the overweighting of small probability
events, which are inconsistent with VNM utility maximization.
To overcome a violation of first-order stochastic dominance in 
prospect theory, cumulative prospect theory (CPT) was
introduced~\cite{tversky1992advances}, which replaces probabilities of outcomes
with their rank-dependent cumulative probability distribution. This change leads to an
overweighting of extreme low-probability outcomes, instead of all
low-probability outcomes. Maximizing CPT utility yields more realistic
predictions of actual human decision-making behavior than maximizing VNM utility.

\subsection{Portfolio optimization}\label{s-port-opt}
Our focus is portfolio optimization, \ie,
choosing a mix of investments in a set of assets.
One approach is based on VNM utility, where the expected value of 
a concave increasing utility function (which is also concave) is
maximized subject to the constraints on the portfolio.
Another approach, introduced by Markowitz~\cite{Markowitz1952},
poses the problem as a bi-criterion optimization problem, with the goal of 
trading off the maximization of expected return with minimization of risk,
taken to be the variance of the portfolio return.
The standard approach is to combine the return and risk, scaled by a
risk aversion factor, into a risk-adjusted return, and maximize this 
concave quadratic objective subject to the constraints. For this reason
Markowitz portfolio optimization is
also referred to as mean-variance (MV) portfolio optimization.
These two approaches are not the same, since the MV 
utility function is not increasing, but they are closely related.  
For example, with a Gaussian asset return
model and exponential utility, VNM portfolio optimization is the same as 
MV portfolio optimization~\cite{merton1969exponential}.
In other cases, MV portfolio optimization was shown to be
approximately optimal for other forms of utility functions~\cite{LevyMarkowitz}.

One advantage of the MV formulation is that the objective can be 
expressed explicitly as a quadratic function, without an expectation over
the asset returns. (The MV objective is the expected value
of a function of the return, but one with a simple analytical expression.)
This enables it to be solved analytically for special cases~\cite{GrinoldKahn},
and very efficiently
using numerical methods for convex optimization when the 
constraints are convex~\cite{cvxbook}.
Leveraging convex optimization,
many extensions were developed, such as the inclusion of transaction and 
holding costs, or multi-period optimization~\cite{BoydKahnMultiPeriod},

VNM portfolio optimization generally uses sample based stochastic 
convex optimization~\cite{shapiro}.  Here we use samples of the 
asset returns, which can be historical or generated from a stochastic model
of asset returns (presumably fit to historical data).
While these methods solve the problem globally, 
they are substantially slower
than MV methods for similar problems.

We now describe these portfolio optimization methods in more detail.
We consider a set of $n$ assets.  A portfolio is characterized by its
asset weights $w=(w_1, \ldots, w_n)$, where $w_i$ is the fraction of the total
portfolio value (assumed to be positive) invested in asset $i$, with $w_i<0$
denoting a short position.
The goal in portfolio optimization is to choose $w$.
The general portfolio optimization problem is 
\BEQ\label{e-port-opt}
\begin{array}{ll}
  \mbox{maximize} & U(w)\\
  \mbox{subject to} & \ones^T w=1, \quad w\in\mathcal{W},
\end{array}  
\EEQ
with variable $w \in \reals^n$.
Here $U:\reals^n \to \reals$ is a utility function,
$\mathcal W \subseteq \reals^n$ is the set of allowable portfolio weights,
and $\ones$ is the vector with all entries one.
The data in this problem are the utility function $U$ and the
portfolio constraint set $\mathcal W$, which we assume is convex.

When $U$ is a concave function the
portfolio optimization problem~\eqref{e-port-opt} is convex,
and so readily solved globally~\cite{cvxbook}.
When $U$ is not concave, the problem
is not convex, and in general difficult to solve globally.  In this case,
we typically resort to heuristic or local methods, which attempt to
solve~\eqref{e-port-opt}, but cannot guarantee that the globally optimal
portfolio is found.

MV portfolio optimization uses the concave quadratic utility function
\[
U^\mathrm{mv}(w)=w^T\mu-\gamma w^T\Sigma w,
\]
the risk-adjusted expected return, where $\mu$ is the expected 
asset return, $\Sigma$ is the covariance matrix of the asset returns, 
and $\gamma >0$ is the risk aversion parameter, used to control
the trade-off of the mean and variance of the portfolio return.
This yields a convex optimization problem that is efficiently solved.

VNM utility maximization uses a utility function of the form
\[
U^\mathrm{vnm}(w) = \Expect u(r^Tw),
\]
where $r$ is the random asset return vector, and $u:\reals \to \reals$ is
a concave increasing utility function. 
Since expectation preserves concavity,
$U$ is a concave function of $w$ and this too leads to a convex portfolio 
optimization problem.  
In a few cases, the expectation can be worked out analytically, but in most
cases one substitutes a sample or empirical average
for the expected value, leading to the approximation
\[
U^\mathrm{vnm}(w) \approx \frac{1}{N} \sum_{i=1}^N u(r_i^Tw),
\]
where $r_1, \ldots, r_N$ are samples of returns.
Maximizing this approximation results in a convex portfolio 
optimization problem.

\subsection{This paper}
We consider portfolio optimization under the 
CPT utility, which we denote $U^\mathrm{cpt}$, defined later in \S\ref{s-cpt}. 
It is well known that CPT utility is not a concave function,
so the problem of choosing portfolio weights so as to maximize it 
is not a convex optimization problem, as VNM utility maximization
and MV portfolio optimization are.
This makes it a challenge to carry out CPT utility 
maximization in practice.

While CPT utility is not concave, we will show that it does have some convexity
structure. Specifically, it is the composition of a convex increasing function
of concave functions for positive returns, and the composition of a concave
increasing function of convex functions for negative returns.
This observation allows us to construct a concave lower bound, or minorant, for
the CPT utility, and leads immediately to a simple algorithm for maximizing it
by repeatedly maximizing the constructed minorant (which is a 
convex problem, and thus readily solved).  This simple 
minorization-maximization (MM) method leads to a local maximum of
the CPT utility~\cite{ccp-lipp}.

Our MM method scales to medium size problems, with perhaps tens of assets
and hundreds of return samples.
For larger problems, we give two other algorithms.
One algorithm uses a simpler optimization of a minorant to the approximation
given by fixing the probability weights that arise in CPT in each step,
again relying on iterations that involve solving convex optimization problems.
As a result, this method can handle complex portfolio constraints, as long as
they are convex.
The second additional algorithm scales to very large problems, but handles
only simple portfolio constraints.
It relies on modern frameworks for 
automatic differentiation and first-order optimization methods.

Open-source Python implementations of all three methods
can be found in the code repository
\url{https://github.com/cvxgrp/cptopt}.

We do not address questions such as whether or when 
one \textit{should} choose a portfolio that maximizes CPT utility.  
We only address the question
of \emph{how} it can be done, algorithmically and computationally. 
We provide methods to solve the CPT portfolio optimization problem, but we emphasize
that the conceptual framework behind the method is more general. Across
computational economics and finance, nonconvex problems frequently arise. The
framework for decomposing the problem into convex and concave parts can be 
extended to such other problems.

\subsection{Previous and related work}

Limited prior work exists on portfolio optimization with CPT utility.
Analytical solutions exist for special cases such as single-period settings with
one risk-free and one risky
asset~\cite{Bernard2010, he2011portfolio, zou2017optimal} or for two-fund
separation under elliptical distributions~\cite{pirvu2012multi}. Extensions
of these special cases to a multi-period setting are considered
in~\cite{shi-multi-period}.
For the general multi-asset cases, heuristics such as particle swarm
simulation~\cite{barro2020cumulative}, or grid search
methods~\cite{hens2014cumulative} are employed, which have been extended to the
multi-period case using dynamic
programming~\cite{degiorgi-dynamic, barberis-loss}. While grid search can
accommodate constraints, the particle swarm method used
in~\cite{barro2020cumulative} cannot, and requires a hyperparameter to be chosen
to turn constraints into penalties.

The evaluation of CPT utility along the mean-variance frontier is a commonly
used heuristic~\cite{LevyLevy2004, hens2014cumulative}.
Some authors (\eg,~\cite{srivastava2022})
use numerical methods to maximize CPT utility on small problems,
do not explicitly mention the numerical solve method,
suggesting the use of generic nonlinear optimizers.
In contrast, we focus on custom methods that exploit the special
structure of the CPT utility maximization problem.

After the initial release of this manuscript, 
Yan et al.~\cite{cpt-admm} proposed a method for optimizing a portfolio
using CPT utility based on the
alternating direction method of multipliers (ADMM) (see, \eg,
\cite{admm}).
This work is closely related to ours, in that they consider 
general multi-asset portfolios, and exploit convexity structure,
although in a different way than we do. It is the only other method we are aware
of that exploits the convexity structure of the CPT utility and can handle constraints.
While we make an approximation about the monotonicity of the weights, they
employ a method which does not globally solve one of their sub-problems in order
to obtain tractable speeds. Carrying out a direct comparison of the methods is not immediately possible

\subsection{Outline}

We start in~\S\ref{s-cpt} by defining CPT utility, fixing our notation.
The CPT utility extends prospect theory (PT) utility, 
described in~\S\ref{s-prospect-util}, by adding a
reweighting function, described in~\S\ref{s-reweighting}.
In~\S\ref{s-cpt-convexity} we explore the convexity structure of CPT utility,
followed by a description of the CPT utility portfolio optimization problem
in~\S\ref{s-cpt-opt}.
In~\S\ref{s-methods} we describe algorithms that can be used to find a
portfolio that maximizes CPT utility.
The first method, presented in \S\ref{s-mm},
is a minorization-maximization method that relies on the convexity structure
described in the previous section.
The second method, described in~\S\ref{s-cc}, uses the convex-concave procedure,
a method for maximizing the sum of a convex and concave function,
and iterates over the probability weights that appear in the CPT utility.
The last method, given in~\S\ref{s-grad},
is a projected gradient type method, which can scale to large problem sizes.
Numerical experiments are presented in~\S\ref{s-examples}, where we evaluate
all methods on a toy problem with three assets,
a medium-sized problem with more assets, and a large-scale problem,
all based on historical asset class data.
We give some conclusions in~\S\ref{s-conclusions}.

\section{Cumulative prospect theory utility}\label{s-cpt}

\subsection{Prospect theory utility}\label{s-prospect-util}

In this section we introduce PT utility, the first
building block of CPT utility. Like VNM utility,
it is monotonically increasing, but PT utility is not concave.
PT utility has an inflection point at the origin, which represents a reference
wealth.
It is concave for positive arguments, \ie,
investors are risk averse for gains, and convex on for negative arguments, \ie,
investors are risk seeking for losses.
Exponential utility functions are commonly used for both the convex and the
concave sections of the PT utility function.
We thus define the positive and negative exponential utilities as
\[
  u_+(x) = 1-\exp(-\gamma_+ x),\qquad u_-(x) = -1+\exp(\gamma_- x),
\]
where $\gamma_+,\gamma_->0$ are parameters. Here and throughout the paper,
functions with a subscript plus sign are applied to gains, and functions
with a subscript minus sign are applied to losses. Combining both functions
yields the exponential prospect theory utility for a single return
\[
u^\mathrm{prosp}(x)=\left\{\begin{array}{ll}
  u_+(x) & \mbox{if }x\geq 0\\
  u_-(x) & \mbox{otherwise}
\end{array}\right.,
\]
which is S-shaped. Prospect theory further accounts for loss aversion, which
requires $\gamma_- > \gamma_+$, \ie, a marginal decrease in wealth would decrease the
utility more than a marginal increase in wealth would increase the utility.

\subsection{Probability reweighting}\label{s-reweighting}

The second building block of CPT utility is a reweighting function that
assigns higher weights to more extreme outcomes. As is common in the CPT literature,
we first define the weighting functions $w(p): [0, 1] \to [0, 1]$.
We take the specific weighting functions
\[
  w_+(p) = \frac{p^{\delta_+}}{(p^{\delta_+}+(1-p)^{\delta_+})^{1/{\delta_+}}}
  ,\qquad
  w_-(p) = \frac{p^{\delta_-}}{(p^{\delta_-}+(1-p)^{\delta_-})^{1/{\delta_-}}},
\]
where $\delta_+,\delta_->0$ are parameters. We now specify the notion of
extreme outcomes. Let $r_1, \ldots, r_N\in \reals^n$ be the empirical
distribution of realized returns on $n$ assets.
Consider a vector of portfolio weights $w\in\reals^n$, with $\ones^Tw = 1$,
where $\ones$ is the vector with all entries one.
The associated portfolio returns are $r_1^Tw, \ldots, r_N^Tw \in \reals$.
Without reweighting, all returns would have equal weight.
Let $N_-$ denote the number of negative returns, and $N_+$ the number of 
nonnegative returns, with $N_-+N_+=N$.
We let $\rho_i$ denote the returns re-ordered or sorted by the portfolio
returns, with index value
$i = 1, \ldots, N$,
\[
w^T\rho_{1}\leq \cdots \leq w^T\rho_{N_-}< 0 \leq w^T\rho_{N_- + 1}\leq
\cdots\leq w^T\rho_{N},
\]
\ie, $w^T\rho_{1}$ is the largest loss and $w^T\rho_{N}$ is the largest gain.
We define the positive and negative decision weights respectively as
\[
  {\pi'_{+,j}} = \left\{
  \begin{array}{ll}
    w_+((N_+ - j +1)/N) - w_+((N_+ - j)/N) & j=1,\ldots,N_+-1 \\
    w_+(1/N) & j=N_+,
  \end{array}  
  \right.
\]
\[
  {\pi'_{-,j}} = \left\{
  \begin{array}{ll}
    w_-((N_--j+1)/N) - w_-((N_--j)/N) & j=1,\ldots,N_--1 \\
    w_-(1/N) & j=N_-.
  \end{array}
  \right.
\]

We would argue that $\pi'_+$ and $\pi'_-$ should be nondecreasing,
\ie, we should put higher weight on more extreme events.
This occurs for most reasonable choices of parameters, but there are choices
for which monotonicity is (slightly) violated.
Thus, we force monotonicity by replacing
${\pi'_{+,j}}$ with $\min(\pi'_+)$ for all $j < \argmin(\pi'_+)$, and likewise
for $\pi'_-$.
We zero-pad $\pi'_+$ and $\pi'_-$ from the left to be length $N$, \ie,
$\pi_+ = (0_{N_-},\pi'_+)$ and $\pi_- = (0_{N_+},\pi'_-)$,
where the subscript on the vector zero denotes its dimension.
We define for a monotone increasing probability vector $\pi$
\[
  f_\pi(x) = \sum_{i=1}^N\pi_ix_{(i)},
\]
which is sometimes called the weighted-ordered-sum or dot-sort function. 
(The notation $x_{(i)}$ means the $i$th smallest element of the vector $x$.)
Then, with
\[
  \phi_+(x) = \max(x,0),\qquad \phi_-(x) = -\min(x,0),
\]
we have the total CPT utility given by
\[
U^\mathrm{cpt}(w) = f_{\pi_+}\left(\phi_+\left(u_+(Rw)\right)\right)-
f_{\pi_-}\left(\phi_-\left(u_-(Rw)\right)\right).
\]

\subsection{Convexity properties}\label{s-cpt-convexity}
In this section we describe some convexity properties of the CPT utility
function. 
PT utility is convex for negative arguments and concave for positive
arguments by definition.
CPT utility, \ie, with reweighting, is a difference of two structured terms
\[
U^\mathrm{cpt}(w) = 
  \underbrace{f_{\pi_+}(\phi_+(}_{\mathrm{convex}}
  \underbrace{u_+(Rw)}_{\mathrm{concave}}))
  -
  \underbrace{f_{\pi_-}(\phi_-(}_{\mathrm{convex}}
  \underbrace{u_-(Rw)}_{\mathrm{convex}})).
\]
The first term is a composition of dot-sort-positive, $f_{\pi}\circ \phi_+$, and
the concave exponential utility for gains, $u_+$.
The dot-sort-positive function is convex, because dot-sort is convex and increasing
for positive weights $\pi$, and $\phi_+$ is convex.
The weighted sum in the CPT utility is consistent with dot-sort whenever
the weights in the dot-sort function are monotone nondecreasing, \ie,
$\pi_1\leq\pi_2\cdots\leq\pi_N$.
Similarly, $f_{\pi}\circ \phi_-$ is convex following the same reasoning,
making $-f_{\pi}\circ \phi_-$ concave, which is in turn composed with the
convex exponential utility for losses, $u_-$.
We note that for each return, only one argument of the difference contributes
to the CPT utility, as $\phi_+(x)\phi_-(x)=0$.
These convexity properties motivate principled algorithmic approaches to
maximizing the CPT utility, which we explore in \S\ref{s-mm} and \S\ref{s-cc}.

\subsection{CPT utility portfolio optimization problem}\label{s-cpt-opt}
The CPT utility portfolio optimization problem is
\BEQ\label{e-cpt-opt}
\begin{array}{ll}
  \mbox{maximize} & U^\mathrm{cpt}(w)\\
  \mbox{subject to} & \ones^{T}w=1,\quad w\in\mathcal{W},
\end{array}
\EEQ
with variable $w$, where $\mathcal{W}$ is the set of
feasible portfolio weights.
It is not a convex optimization problem, so we will seek approximate
solution methods.

We mention some simple methods for solving or approximately solving
the CPT utility portfolio optimization problem~\eqref{e-cpt-opt}.
If the number of assets is very small (say, 3 or 4), we can solve it
by brute force computation,
by evaluating $U^\mathrm{cpt}$ over a fine grid of values.

A reasonable heuristic for approximately solving the CPT utility portfolio
optimization problem, motivated by~\cite{LevyLevy2004}, leverages our ability to
efficiently solve the MV portfolio optimization problem.
We find the so-called efficient frontier, by solving the MV problem
for a number of different values of the risk aversion parameter $\gamma$.
(This gives the MV efficient frontier.)
We evaluate the CPT utility of each of these portfolios, and choose the one with the
largest value.  While this does not in general solve the
problem~\eqref{e-cpt-opt}, it often produces a very good, \ie,
nearly optimal, portfolio.  It can be used as an initial guess for
the iterative methods described below.
We refer to this method as the MV heuristic for CPT maximization.

\section{Optimization methods}\label{s-methods}
\subsection{Minorization-maximization method}\label{s-mm}
The CPT utility has the composition form
\[
U^\text{cpt}(w) = (f_{\pi_+}\circ \phi_+)(u_+(Rw))
-(f_{\pi_-}\circ \phi_-)(u_-(Rw)).
\]
We denote a general linearization of a convex (concave) function $h(w)$ at the
point $\hat{w}$ as
\[
  \widehat{h}(w, \hat{w}) = h(\hat{w}) + g^T (w - \hat{w}),
\]
where $g$ is a subgradient (supergradient) of the
function $h$.
As all linearizations that follow occur at $\hat{w}$, we suppress the second
argument.

At $\hat{w}$, we create a concave approximation of the first term of
$U^\text{cpt}(w)$ by linearizing $(f_{\pi_+}\circ \phi_+)$.
We approximate the second term in the difference by linearizing the inner
convex utility $u_-$.
To linearize dot-sort-positive, we observe that a subgradient is
given by the vector $g_x$ with entries
\[
    g_{x,i} = \left\{\begin{array}{ll}
    0 & \mbox{if } x_i<0\\
    {\pi_+}_{\sigma_x(i)} & \mbox{otherwise,}
\end{array}\right.
\]
where $\sigma_x$ is the permutation which maps $i$ to the rank of $x_i$ in $x$.
The minorant at $\hat{w}$ is therefore given by
\[
\tilde{U}^\text{cpt}(w) = (\widehat{f_{\pi_+}\circ \phi_+})(u_+(Rw))
-(f_{\pi_-}\circ \phi_-)(\widehat{u_-}(Rw)).
\]
The minorization-maximization (MM) algorithm (also called the majorization-minimization
algorithm when solving a minimization problem) simply iterates between
creating the minorant at the current iterate and then maximizing it to find 
the next iterate~\cite{MM}.  Our minorant is concave, so maximizing it is
efficient.
Here $\mathcal{W}$ can be any DCP convex constraint set, since each iteration
requires solving a general convex optimization problem.

\begin{algdesc}{Minorization-maximization method} \label{alg-mm}

\textbf{given} $w^0$, let $\hat{w} := w_0$.
\\\textbf{repeat:}
\begin{enumerate}
 \item Let $w^\text{next}$ be a maximizer of $\tilde{U}^\text{cpt}(w)$,
 subject to $w\in\mathcal{W}$.
 \item {\bf break if} $w^\text{next}=\hat{w}$.
 \item Update $\hat w=w^\text{next}$.
\end{enumerate}
{\bf return} $\hat w$.
\end{algdesc}

\subsection{Iterated convex-concave method}\label{s-cc}
Though the CPT portfolio optimization objective is non-convex, we know the
curvature and sign properties of the component functions which are composed to
form the utility. In particular, PT utility is convex on the negative reals,
and concave on the nonnegative reals.
Thus, it is amenable to optimization via the convex-concave
procedure (CCP)~\cite{ccp-lipp, dccp, ccp-yuille,ccp-lanckriet}.
The convex-concave procedure for maximization iteratively linearizes the second
term in the sum of a concave and a convex function, and maximizes this surrogate
objective.
While the PT utility has this clear curvature, the CPT utility does not, due
to reweighting.
Our heuristic approach is to fix the probability weights in each iteration,
and then solve the fixed weight CPT utility optimization problem with the
convex-concave procedure.
Once the weights have been fixed, we can write the PT utility as a concave
function
\[
f^{\text{ccv}}(x) = \left\{\begin{array}{ll}
  1-\exp(-\gamma_+ x)  & \mbox{if } x\geq 0\\
  \gamma_-x & \mbox{otherwise},
\end{array}\right.
\]
plus a convex function,
\[
  f^{\text{cvx}}(x)=\inf_{z\leq 0,\; z\leq x} \left( -1+\exp(\gamma_-z)-\gamma_-z\right).
\] 
Here ``cvx'' and ``ccv'' denote convex and concave, respectively. (See \S\ref{s-append-dcp} for a derivation of these functions in disciplined
convex programming (DCP) form.) 
Unlike the MM algorithm in \S\ref{s-mm} which maximizes a global lower
bound, this approximation is only local, so we include a trust region
constraint, which we omit from the algorithm description for brevity. Note that
as before, $\mathcal{W}$ can be any DCP convex constraint set.

\begin{algdesc}{Convex-concave procedure}\label{alg-cc}

\textbf{given} $w_0$, let $\hat{w} := w_0$.
\\\textbf{repeat:}
\begin{enumerate}
  \item Let $\pi$ be the concatenation of the reversed vector $\pi'_-$ followed by $\pi'_+$,
  where $\pi'_-$ and $\pi'_+$ are the decision weights associated with $\hat{w}$
  (see \S\ref{s-cpt}).
  \item Let $L^{\text{cvx}}_i$ be the linearization of $f^{\text{cvx}}$ at $(R\hat{w})_{(i)}$
  \item Let $w^\text{next}$ be a maximizer of $\sum_i \pi_i \left(f^{\text{ccv}}\left(w^T\rho_i\right)
  -L_i^{\text{cvx}}\left(w^T\rho_i\right)\right)$,
  subject to $w\in\mathcal{W}$, where $\rho_{i}$ is the row of $R$ associated with $(R\hat{w})_{(i)}$.
 \item {\bf break if} $w^\text{next}=\hat{w}$.
 \item Update $\hat w=w^\text{next}$.
\end{enumerate}
{\bf return} $\hat w$.
\end{algdesc}

\subsection{Projected gradient ascent}\label{s-grad}
We first consider maximizing the CPT utility using gradient ascent (GA).
While the CPT utility is not differentiable everywhere, we can use an automatic
differentiation package such as PyTorch~\cite{torch2019}
to specify the computation chain for the
problem and automatically compute the gradient at points where the utility is
differentiable, and a reasonable surrogate for the gradient
(such as a subgradient for convex functions) at points
where it is not.
Such libraries are extremely fast and optimized for use on GPUs.
We can then perform gradient ascent, together with a method to enforce the
portfolio constraints.
Projected gradient ascent consists of the iterations
\[
  w^{k+1}=\Pi\left(w^k+\eta^k\nabla f(w^k)\right),
\]
where $k$ denotes iteration, $\eta^k>0$ is a stepsize and 
$\Pi$ is $\ell_2$ or Euclidean
projection onto the constraint set $\mathcal{W}$, \ie,
$\Pi(w)=\argmin_{w'\in\mathcal{W}} \|w'-w\|_2$.

Using these computation frameworks requires the projection to
be expressed as a simple computation chain,
which can be done in simple cases such as a long-only portfolio, \ie,
$\mathcal W=\reals_+^n$.
Another option to handle long-only portfolio constraints is to parametrize
nonnegative portfolio weights via a multinomial logistic map,
\[
w_i=\frac{\exp x_i}{\sum_j\exp x_j}, \quad i=1, \ldots, n,
\]
where $x$ is an unconstrained variable.

\clearpage
\section{Numerical examples}\label{s-examples}
To evaluate the efficiency and performance of the proposed methods,
we compare them in a series of numerical experiments with increasing data size.
We first compile a data set consisting of $N=600$ monthly returns,
covering the 50-year period from 07--1972 to 06--2022. The $n=14$ assets consist
of equities and fixed income securities from different regions, as well as
commodities, as displayed in table~\ref{t-asset-classes}. 
The data was obtained from Global Financial Data (GFD)~\cite{cfd-data}.
The indices are total return indices, \ie, they include dividends and interest
payments. The GFD index methodology extends the index history back in time
by combining multiple single indices where necessary.
We provide the GFD symbol for each asset class for reference.

\begin{table}
  \begin{minipage}{\textwidth}
    \centering
    \begin{tabular}{|c|c|c|}
      \hline
      Asset class      & Region           & GFD symbol \\ \hline
      Equity           & US               & SPXTRD \\
      Equity           & Europe           & STOXXER \\
      Equity           & Japan            & TOPXDVD \\
      Equity           & Emerging markets & TRGFDEM \\
      Government Bonds & US               & TRUSG10M \\
      Corporate Bonds  & US               & TRCCRBD \\
      Government Bonds & Europe           & TREUROGM \\
      Government Bonds & Japan            & TRJPNGVM \\
      Bills            & US               & TRUSABIM \\
      Bills            & Europe           & TREUROBM \\
      Bills            & Japan            & TRJPNBIM \\
      Commodities      & Global           & TRUSACOM \\
      Gold             & Global           & XAU\_BD \\
      Silver           & Global           & XAG\_HD \\
      \hline
    \end{tabular}
    \caption{Asset classes and regions in the data set.}
    \label{t-asset-classes}
  \end{minipage}
\end{table}

\subsection{Toy example}
Our first small example uses $n=3$ assets: US stocks, 10-year US Treasury bonds,
and 3-month US Treasury T-bills.
We choose the CPT function with parameters
\[
  \gamma_+ = 8.4, \quad \gamma_- = 11.4,\qquad
  \delta_+ = 0.77,\quad \delta_- = 0.79,
\] 
which are reasonable, and at the same time exhibit clear non-convexity and
even multimodality of the CPT utility.
(Many other reasonable choices of the parameters lead to unimodal CPT utility,
which makes the portfolio optimization problems easy; our goal is to 
evaluate the methods on more challenging problem instances.)

\begin{figure}
  \begin{center}
    \includegraphics[width=.7\textwidth]{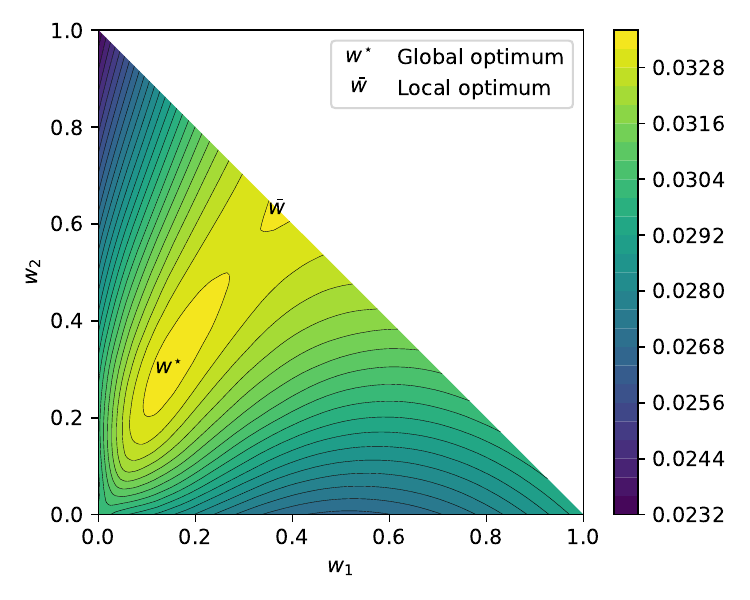}
    \caption{CPT utility surface for a long-only portfolios of
    stocks ($w_1$), bonds ($w_2$), and T-bills ($w_3=1-w_1-w_2$).}
    \label{f-utility-surface}
  \end{center}
\end{figure}

Figure~\ref{f-utility-surface} gives a plot of this utility function for
$\mathcal{W}=\reals^2_+$, \ie, long-only portfolios.  The horizontal axis
is $w_1$, the fraction invested in stocks; the vertical axis is $w_2$, 
the fraction invested in bonds. The fraction invested in T-bills is
$w_3 = 1-w_1-w_2$. Thus, the point $(0,0)$ represents a portfolio
fully invested in T-bills. Any portfolio on the diagonal connecting
$(1,0)$ and $(0,1)$ represents portfolios invested in a convex combination of
only stock and bonds.  Since there are only two portfolio weights to optimize over,
we can find the
global maximum using brute-force evaluation of the utility over a fine grid.
The global maximum is attained at $w^{\star}=(0.14, 0.3)$,
yielding $U^\mathrm{cpt}(w^{\star})=0.0334$. 
In addition, there is a local optimum with
slightly lower utility at $\bar{w}=(0.37, 0.63)$, which yields
$U^\mathrm{cpt}(\bar{w})=0.0332$.

\begin{figure}
  \begin{center}
    \includegraphics[width=1\textwidth]{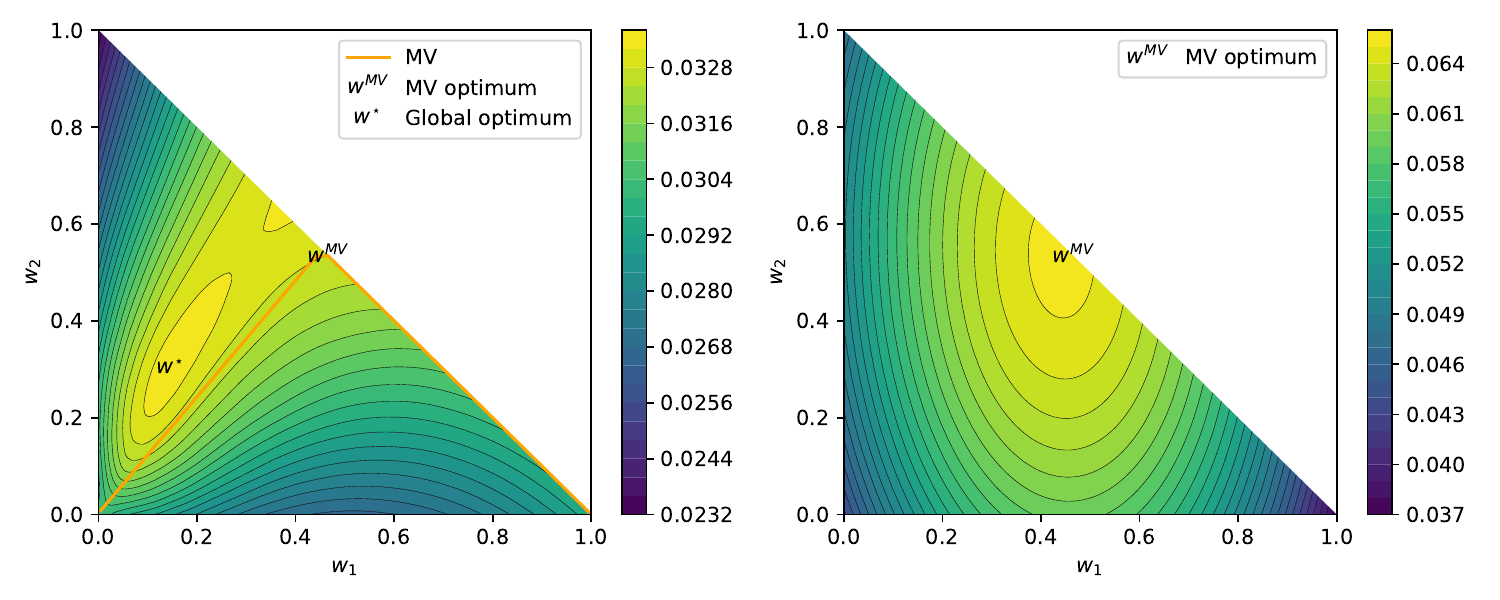}
    \caption{(a) Maximizing $U^{\mathrm{cpt}}$ along the MV frontier, resulting
    in $w^{\mathrm{mv}}$. (b) Utility surface of $U^{\mathrm{mv}}$ for the
    choice of $\lambda$ that results in $w^{\mathrm{mv}}$.}
    \label{f-markowitz}
  \end{center}
\end{figure}

\paragraph{MV frontier.}
A simple heuristic is to evaluate $U^{\mathrm{cpt}}$ on portfolios along the
mean-variance efficient MV frontier and choosing the maximizing portfolio among
them. Based on the sample mean and covariance of the returns, we first find
the return-maximizing and risk-minimizing portfolios, and then sample 100 points
that are equidistant in volatility space along the efficient MV frontier.
Figure~\ref{f-markowitz}~(a) shows the efficient MV frontier, and the portfolio with the
highest CPT utility along it, $w^{\mathrm{mv}}$, associated with risk aversion
parameter $\gamma =3.2$. It achieves CPT utility
of $U^\mathrm{cpt}(w^{\mathrm{mv}}) = 0.0328$.
Figure~\ref{f-markowitz}~(b) shows $U^{\mathrm{mv}}$ for the choice $\gamma =3.2$.
It should be noted that the MV frontier is independent of the parameter
choice of the CPT utility function, and in general the MV optimum can be far
away from a local optimum of CPT\@.

\begin{figure}
  \begin{center}
    \includegraphics[width=0.5\textwidth]{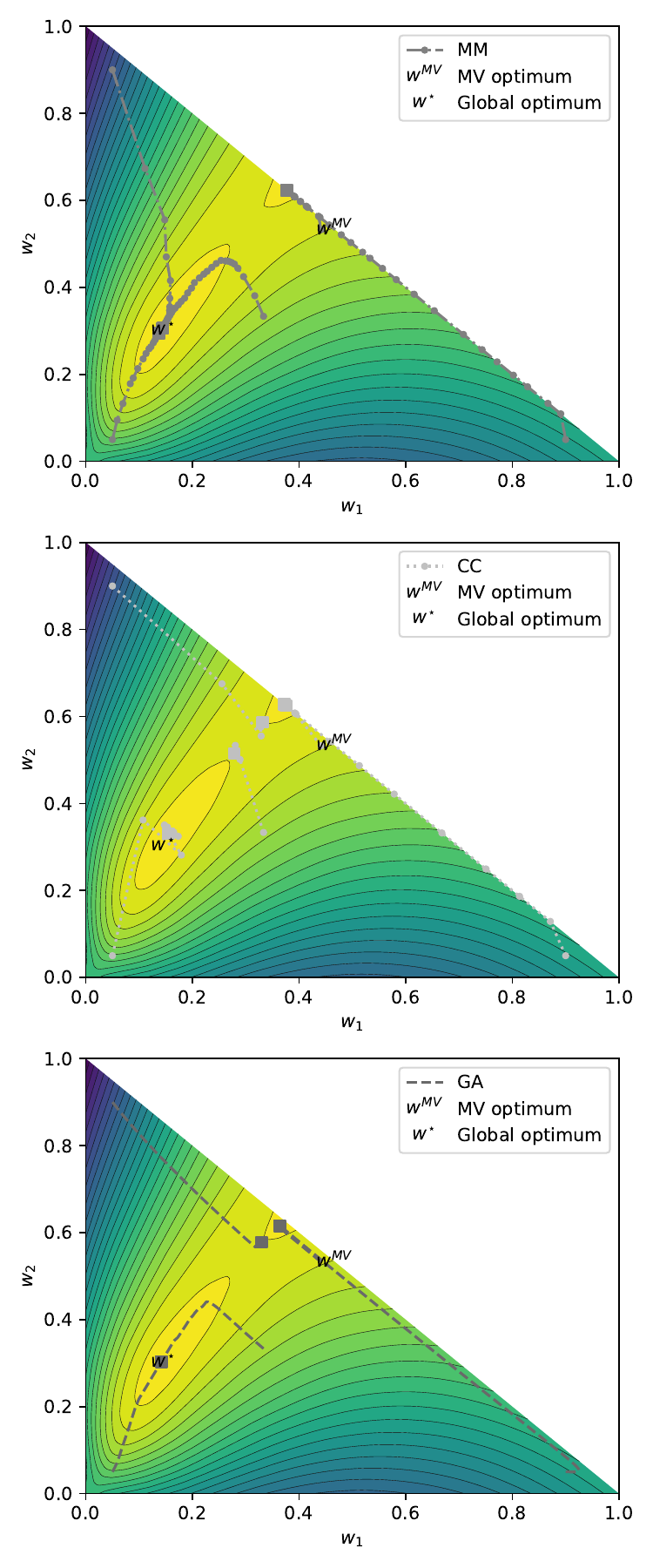}
    \caption{Convergence from different staring points for the MM (a), CC (b),
      and GA (c) methods.}
    \label{f-panel}
  \end{center}
\end{figure}

\paragraph{Iterative methods.}
As all remaining methods depend on initialization, we compare the convergence
from equal weights, three points close to a full investment in each single
asset, as well as the MV optimum in figure~\ref{f-panel}.
The MM algorithm terminates at a local maximum from all starting points within
fewer than 30 iterations.
Likewise, CC converges to a local optimum or a point where the numerical
stopping criterion is reached in all cases within at most 11 iterations,
albeit on a visually more erratic path.
Lastly, the GA method also converges to a point close to a local optimum in all
cases.
Thus, all iterative methods appear to perform equally well on the toy example.

\begin{figure}
  \begin{center}
    \includegraphics[width=0.8\textwidth]{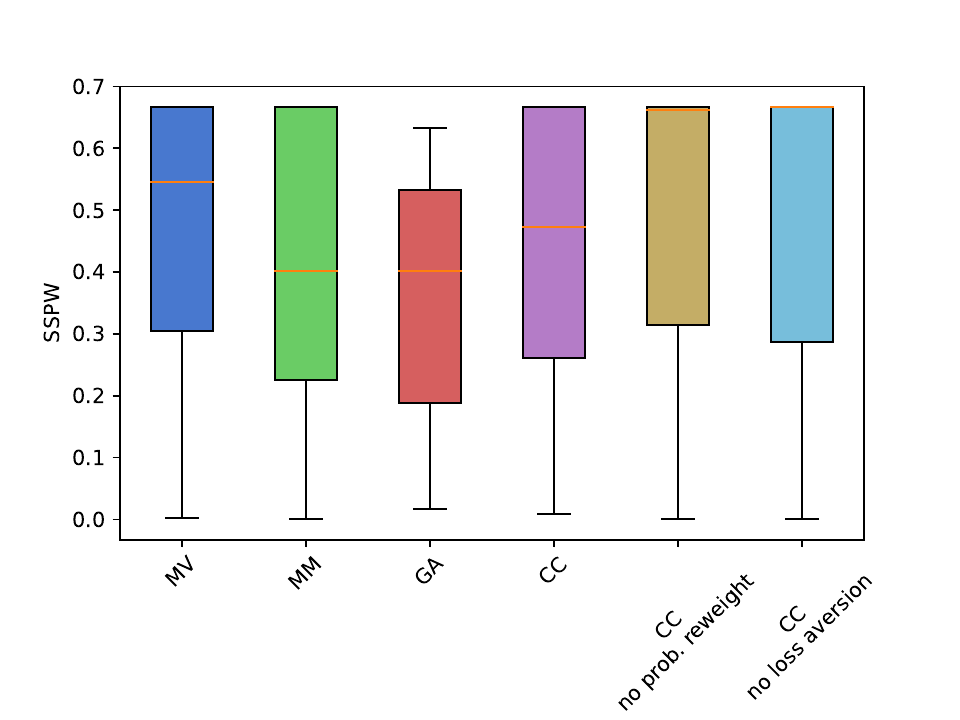}
    \caption{Comparison of the sum of squared portfolio weights across
      methods.}
    \label{f-sspw}
  \end{center}
\end{figure}

\paragraph{Diversification.}
To see the effect of the CPT utility on diversification, as well as to understand
how the different methods alter the portfolio weights, we run a backtest using
a sliding window of 100 observations along our previously described data set of
600 monthly returns.
Following~\cite{spww} and \cite{cpt-admm}, we
compute the sum of squared portfolio weights (SSPW), $\|w - \frac{1}{n}\|_2^2$, as a
measure of diversification. We compare the MM, CC, and GA methods to the MV 
heuristic in figure \ref{f-sspw}. We observe that the MM and CC methods all
result in a lower SSPW than the MV heuristic, \ie, the portfolios are more
diversified. In addition, we explore how the assumptions of 
probability reweighting and loss aversion that are inherent to the CPT utility
affect the portfolio weights. For this, we change set $\gamma_+ = \gamma_- = 11.4$
in the no loss aversion setting, and $\delta_+ = \delta_- = 1$ in the
no probability reweighting setting. We find that both cases result in a 
much higher median SSPW, indicating that the portfolio weights are less 
diversified.

\subsection{Multi-asset example}\label{s-multi-asset-example}
\begin{figure}
  \begin{center}
    \includegraphics[width=1\textwidth]{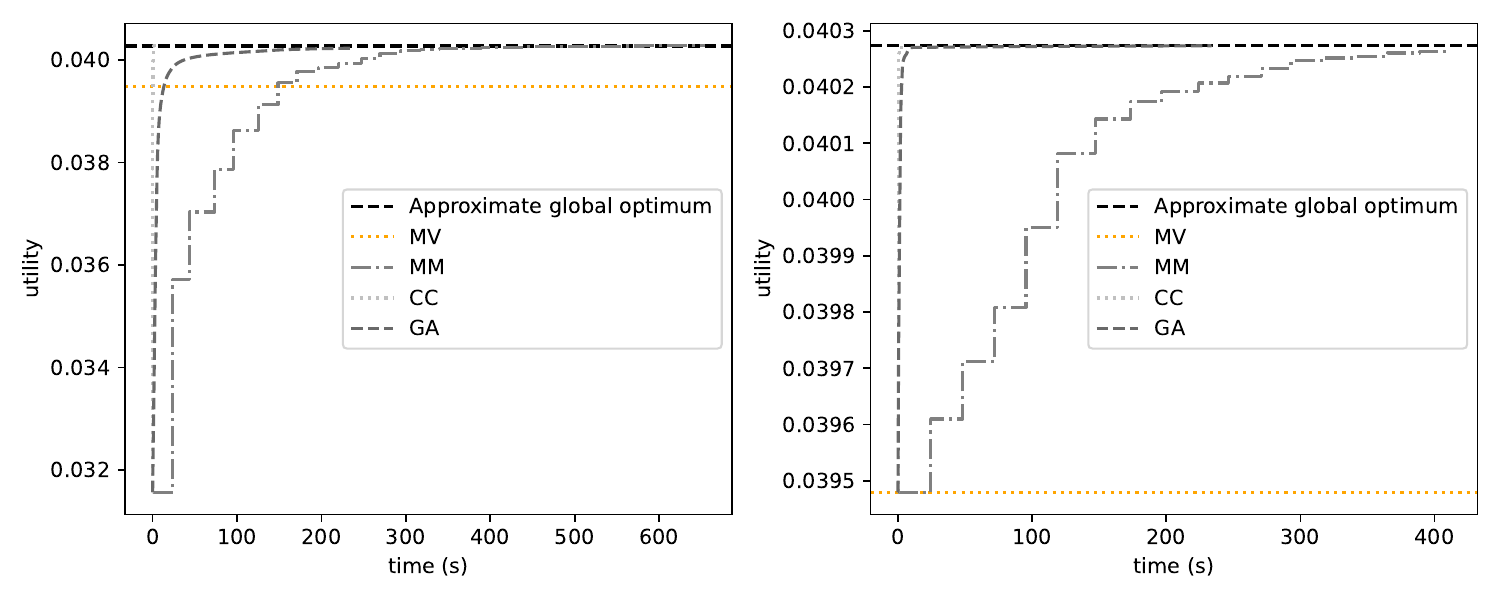}
    \caption{Comparison of wall-time across methods for the multi-asset example,
      started from (a)~equal-weight portfolio and (b)~the best MV portfolio.}
    \label{f-wall-time-medium-panel}
  \end{center}
\end{figure}

We now extend our example to use all $n=14$ assets, comparing the
achieved utilities, as well as the required computation time.
While comparing the absolute wall-times across different implementations can
only approximate the computational efficiency of the algorithms, it is relevant
to the practicality of the presented methods. The best portfolio on the efficient
MV frontier attains a utility of 0.0395 in only 0.6 seconds.
Starting all iterative methods from the equal-weight portfolio, CC terminates
first, yielding a utility of 0.0403 in 4.1 seconds. MM also attains the same
utility, but it takes substantially longer, terminating after about 650 seconds.
GA also results in approximately the same utility, being slower than CC,
but still dominating MM\@.
When optimizing a single portfolio, we find that GA converges faster using the
CPU. However, simultaneous optimization of many portfolios, which is naturally
handled by this method, scales better when using a GPU. We use this observation
to simultaneously optimize from 10,000 starting points. These weights are
sampled from a symmetric Dirichlet distribution with concentration parameter
$\alpha = 1$, \ie, $w_0 \sim \text{Dir}(\ones^n)$, which is equivalent to
a uniform distribution over the open standard $(n-1)$-simplex.
The best resulting utility is denoted as the approximate global optimum, and
is not higher than the utilities achieved by all iterative methods when started
from equal weights. All methods, as well as the approximate global optimum, are
visualized in figure~\ref{f-wall-time-medium-panel}~(a).
In practice, it would likely be beneficial to leverage the
fast computation of the MV portfolio as a starting point for the iterative
methods. Indeed, as shown in figure~\ref{f-wall-time-medium-panel}~(b), this
reduces the time to convergence substantially for all methods, with the GA
method now converging in approximately ten seconds. MM is also faster, but
still takes about 400 seconds to converge.

Investigating the sensitivity to the starting point, we sample 30 starting
weights and compare the attained utilities. We find that MM and CC converge to
a utility of 0.0403 in all cases. GA, however, displays a higher variance.
Its best utility is close to the values obtained by MM and CC\@.
The median utility is 0.0401, which is higher than MV at 0.0395. The worst case
utility is 0.0386, which is worse than all other methods.

\subsection{Scaling to many return samples}
\begin{figure}
  \begin{center}
    \includegraphics[width=.7\textwidth]{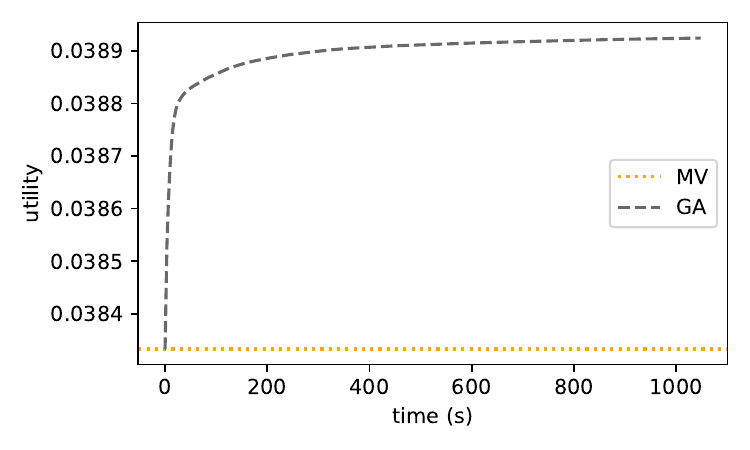}
    \caption{Convergence of the GA method for $N=6,000$ starting from
      $w^{\mathrm{mv}}$.}
    \label{f-large-scale-ga}
  \end{center}
\end{figure}

To investigate the scalability of the methods to more return samples, we extend
the 600 observations of our data set with synthetic returns. For this, we sample
from a Gaussian mixture model with three components that was fit to the return data.
We find that GA scales best, handling data sets of hundreds of thousands
of observations. For such large data sets, MM and CC did not converge in a
reasonable time. Further, as the GA implementation naturally handles optimizing
multiple starting points simultaneously, the problem of high variance
observed in \S\ref{s-multi-asset-example} is alleviated.
Figure~\ref{f-large-scale-ga} shows an example where we extend the original
data set by a factor of ten, \ie, we consider the case where $N=6,000$.
We choose $w^{\mathrm{mv}}$ as a principled starting point for the GA method.
We observe that GA does improve over its starting point, however, looking at
the axis scale reveals that the improvement over MV is marginal.
We observe that the numerical value of the highest utility is different
compared to the original data set in \S\ref{s-multi-asset-example}.
This is expected, however, as the observed data only makes up 10\% of the
extended data set, and the data generating process for the synthetic returns
is only an approximation of the true data generating process, which may not be
fully described by any single distribution.

\section{Conclusions}\label{s-conclusions}
While the CPT utility is nonconvex and can even be
multimodal, we identify some simple convexity properties.
Specifically, the CPT utility is a difference of two structured functions,
with the first term given by a composition of a convex function with concave
arguments and the second term given by a composition of a convex function with
convex arguments.
This structure allows us to construct locally tangent concave minorants,
which we use to develop a minorization-maximization algorithm to maximize the
CPT utility numerically.  While the analysis was restricted to the CPT utility,
we believe that it motivates similar analyses for other nonconvex objectives
commonly used in finance and economics.
We provide several practical methods to maximize the CPT utility, including one
massively scalable method, and two methods which can easily handle arbitrary
convex portfolio constraints. To the best of our knowledge, previous work
on maximizing CPT utility considered only simple analytical cases or small
problem instances with generic nonlinear optimizers.

As a practical matter, for small problems with arbitrary convex constraints,
the MM method has shown smooth convergence and is thus the recommended default
method.
If this method is too slow, but the portfolio constraints are complex,
the CC method should be used instead.
For large problems with simple constraints, the GA method appears to be the best
choice. As there is low scaling overhead, one should optimize many portfolios
simultaneously, including the MV optimal portfolio, an equal weight
portfolio, as well as randomly sampled starting points. As all methods are
readily available in the accompanying code, it is easy to experiment for the
given use case.

Lastly, it is worth noting that the simple method of approximately maximizing
the CPT utility by restricting the feasible set to the MV frontier seems
to closely approximate the optimal CPT utility in many problem instances.

\section*{Declarations}
\paragraph{Funding.} P. Schiele is supported by a fellowship within the IFI program of the German
Academic Exchange Service (DAAD).
This research was partially supported by ACCESS (AI Chip Center for Emerging
Smart Systems), sponsored by InnoHK funding, Hong Kong SAR,
and by ONR grant N00014-22-1-2121.

\paragraph{Competing interests.} The authors have no competing interests to
declare that are relevant to the content of this article. 

\section*{Data availability}
The code is available at \url{https://github.com/cvxgrp/cptopt}.
The dataset analyzed during the current study are not publicly available due
to licensing reasons but are available from the corresponding author on
reasonable request.

\clearpage
\bibliography{cpt.bib}

\clearpage
\appendix
\section{DCP form of CCP objective}\label{s-append-dcp}

\subsection{DCP form of $f^{\text{ccv}}$}
To obtain the piecewise definition of $f^{\text{ccv}}$, we split up its
argument into a positive and negative part,
\[
  x = x^+ + x^-, \qquad x^+ \geq 0, \quad x^- \leq 0.
\]
Now, we observe that
\[
  f^{\text{ccv}}(x)= 1-\exp(-\gamma_+ x^+) + \gamma_-x^-
\]
is concave, because
$\gamma_- > \gamma_+ \geq \frac{\partial (1-\exp(-\gamma_+ x))}{\partial x}$
for $x\geq0$. We implement $f^{\text{ccv}}$ in DCP form via its hypograph,
\[
    \{(x,t)\mid f(x) \geq t \} =
    \{(x,t)\mid \exists~x^+ \geq 0,~x^-\leq0,~ x = x^+ + x^-,~1-\exp(-\gamma_+ x^+) + \gamma_-x^- \geq t\},
\]
which in practice means that to add this function to an optimization problem,
we introduce new variables $t$, $x^+$, and $x^-$,
replace $f^{\text{ccv}}$ with $t$, and add the constraints
\[
    t\leq 1-\exp(-\gamma_+ x^+) + \gamma_-x^-,
    \qquad x = x^+ + x^-,
    \qquad x^+ \geq 0,
    \qquad x^- \leq 0.
\]

\subsection{DCP form of $f^{\text{cvx}}$}

To see that $f^{\text{cvx}}$ is convex, we can equivalently represent it as a
partial minimization of the convex function (of $z$ and $x$ jointly)
\[
  -1+\exp(\gamma_-z)-\gamma_-z + \ones \{ z\leq x \}
\]
over the convex set $\{(z,x) \mid z \leq 0 \}$. The function can be
used in DCP frameworks that provide the indicator function and partial
minimization.
Alternatively, the indicator function can be omitted when adding the
explicit constraints
\[
    z\leq x,
    \qquad z\leq 0.
\]

\newpage
\section{Code snippets}\label{s-code-example}

\begin{lstlisting}[language=Python]
from scipy.stats import multivariate_normal as normal

from cptopt.optimizer import *
from cptopt.utility import CPTUtility

# Generate returns
corr = np.array([
    [1, -.2, -.4],
    [-.2, 1, .5],
    [-.4, .5, 1]
])
sd = np.array([.01, .1, .2])
Sigma = np.diag(sd) @ corr @ np.diag(sd)

np.random.seed(0)
r = normal.rvs([.03, .1, .193], Sigma, size=100)

# Define utility function
utility = CPTUtility(
    gamma_pos=8.4, gamma_neg=11.4,
    delta_pos=.77, delta_neg=.79
)

initial_weights = np.array([1/3, 1/3, 1/3])

# Optimize
mv = MeanVarianceFrontierOptimizer(utility)
mv.optimize(r, verbose=True)

mm = MinorizationMaximizationOptimizer(utility)
mm.optimize(r, initial_weights=initial_weights, verbose=True)

cc = ConvexConcaveOptimizer(utility)
cc.optimize(r, initial_weights=initial_weights, verbose=True)

ga = GradientOptimizer(utility)
ga.optimize(r, initial_weights=initial_weights, verbose=True)
\end{lstlisting}

\end{document}

%% file: defs.tex
\newcommand{\ones}{\mathbf 1}
\newcommand{\reals}{{\mbox{\bf R}}}

\newcommand{\Expect}{\mathop{\bf E{}}}

\newcommand{\argmin}{\mathop{\rm argmin}}

\newcommand{\eg}{{\it e.g.}}
\newcommand{\ie}{{\it i.e.}}

\newcommand{\BEAS}{\begin{eqnarray*}}
\newcommand{\EEAS}{\end{eqnarray*}}
\newcommand{\BEA}{\begin{eqnarray}}
\newcommand{\EEA}{\end{eqnarray}}
\newcommand{\BEQ}{\begin{equation}}
\newcommand{\EEQ}{\end{equation}}
\newcommand{\BIT}{\begin{itemize}}
\newcommand{\EIT}{\end{itemize}}

\newcounter{algorithmctr}
\renewcommand{\thealgorithmctr}{\arabic{algorithmctr}}
\newenvironment{algdesc}%
   {\mbox{}\\*[\parskip]\begin{minipage}{\linewidth}%
       \refstepcounter{algorithmctr}\begin{list}{}{%
       \setlength{\rightmargin}{0\linewidth}%
       \setlength{\leftmargin}{.05\linewidth}}%
       \rmfamily\small
       \item[]{\setlength{\parskip}{0ex}\hrulefill\par%
        \nopagebreak{\bfseries\textsf{Algorithm \thealgorithmctr~}}}}%
   {{\setlength{\parskip}{-1ex}\nopagebreak\par\hrulefill\\*[2ex]\par}%
   \end{list}\end{minipage}}